\documentclass[a4paper,10pt]{article}
\usepackage{psfrag,epsfig,color}
\usepackage{amsmath,amssymb,amsthm}
\usepackage{hyperref}
\usepackage[usenames,dvipsnames]{xcolor}
\usepackage{tcolorbox}

\def \Z{{\mathbb Z}}
\def \R{{\mathbb R}}
\def \C{{\mathbb C}}

\def \H{{\mathbb H}}
\def \Q{{\mathbb Q}}

\newcommand{\Canakci}{\c{C}anak\c{c}\i }
\newcommand{\Ilke}{\.{I}lke }
\newcommand{\I}{\.{I}}

\theoremstyle{definition}

\newtheorem{ex}{Example}
\newtheorem{rem}{Remark}

\begin{document}

\begin{center}
{  \huge \bf 
  Ptolemy Relation and Friends
}

\vspace{10pt}
{\Large  Anna Felikson}

\vspace{10pt}
{\it To Andrei Zelevinsky who would have turned 70 now}
\end{center}

\section{Prelude: Ptolemy's Theorem}
We start with a theorem known since Ptolemy (Claudius Ptolemaeus, 2nd century AD), who used it to create his table of chords with the aim of applying it in astronomy. The statement reads as follows:

%\begin{theorem}[Ptolemy]
\begin{center}
{\it Given a cyclic Euclidean quadrilateral $ABCD$, one has $AB\cdot CD+BC\cdot AD=AC\cdot BD$.}
\end{center}
  % \end{theorem}

\noindent
When $ABCD$ is a rectangle, Ptolemy's theorem turns into the Pythagoras  theorem.

\begin{figure}[!h]
  \begin{center}
    \includegraphics[width=0.99\linewidth]{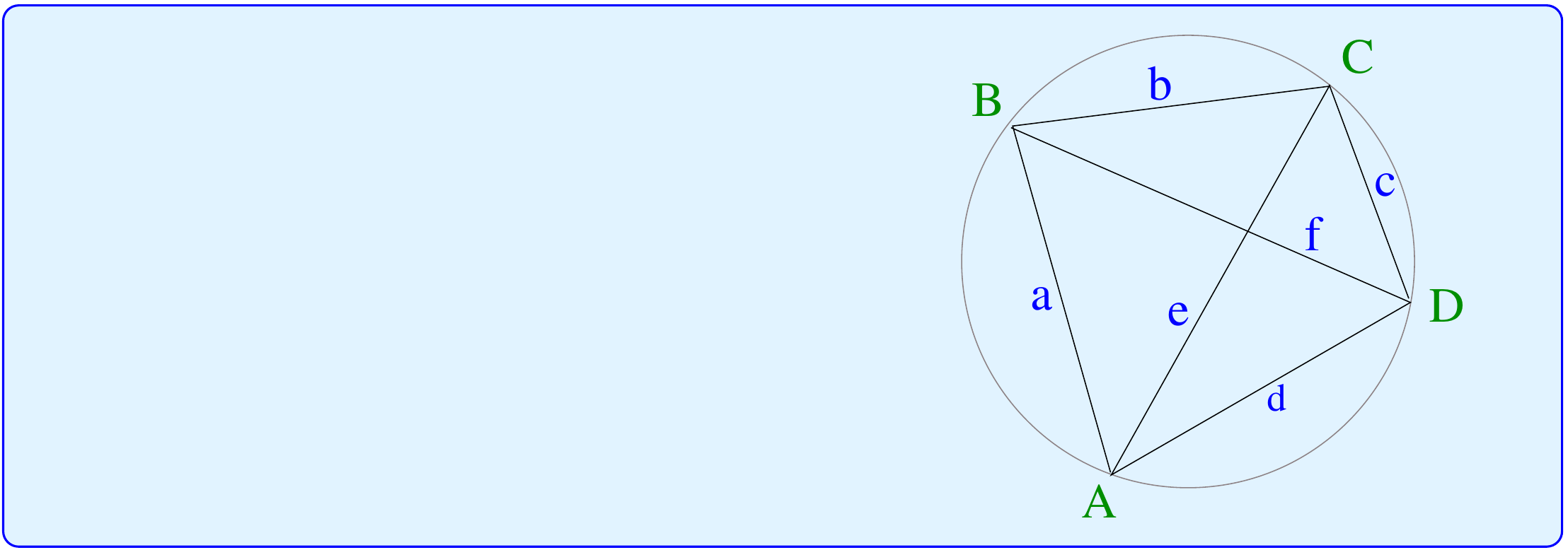}
    \put(-405,15){ \includegraphics[width=0.3\linewidth]{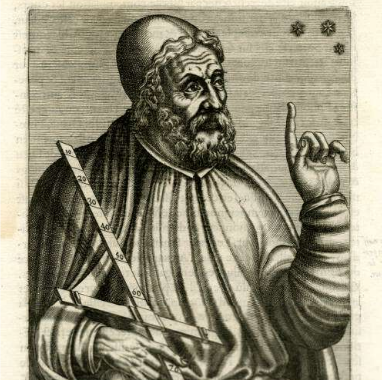}}
    \put(-210,10){\bf \color{blue} ef=ac+bd}
 \caption{Ptolemy (left) and Ptolemy's theorem (right). {\small The \href{https://www.britishmuseum.org/collection/object/P_1879-1213-135}{image} in the left: Illustration by unknown author to ``Les Vrais Pourtraits et vies des Hommes Illustres'', by Andre Thevet, 1584. © The Trustees of the British Museum.} }
\label{perp}
\end{center}
\end{figure}
%Line engraving, French, 1584
%Woodcut from Les Vrais Pourtraits et vies des Hommes Illustres, Andre Thevet, 1584
%https://www.britishmuseum.org/terms-use/copyright-and-permissions
%https://www.britishmuseum.org/collection/object/P_1879-1213-135

In  recent decades,  identities similar to the one in the Ptolemy's theorem started to   pop up in many  fields in connection to the notion of cluster algebras introduced and studied since  2000 by Fomin and Zelevinsky~\cite{FZ1,FZ2}. In this brief note we will try to describe several animals from this big and rich zoo.

\section{Pl\"ucker relations}
Consider the Grassmannian $Gr_{2,n}=\{V \ | \ V\subset \R^n, \ \dim V=2 \}$ of 2-dimensional subspaces in the real $n$-dimensional space. A 2-dimensional subspace $V\subset \R^n$ can be described by two $n$-dimensional vectors $(a_{i1},\dots,a_{in})$, $i=1,2$  spanning the subspace, i.e. by a $2\times n$ matrix  of rank 2.
Denote
$p_{ij}=det\begin{pmatrix}a_{1i} &a_{1j}\\a_{2i} & a_{2j}  \end{pmatrix}  $, the Pl\"ucker coordinates on $Gr_{2,n}$.
The set of these determinants (considered up to simultaneous scaling) completely defines the subspace $V$, so it provides  coordinates on  $Gr_{2,n}$.
It is easy to show that $p_{ij}$ satisfy the Pl\"ucker relation $p_{ik}p_{jl}=p_{ij}p_{kl}+p_{jk}p_{ki}$, for $i<j<k<l$.

\begin{figure}[!h]
  \begin{center}
 \includegraphics[width=0.99\linewidth]{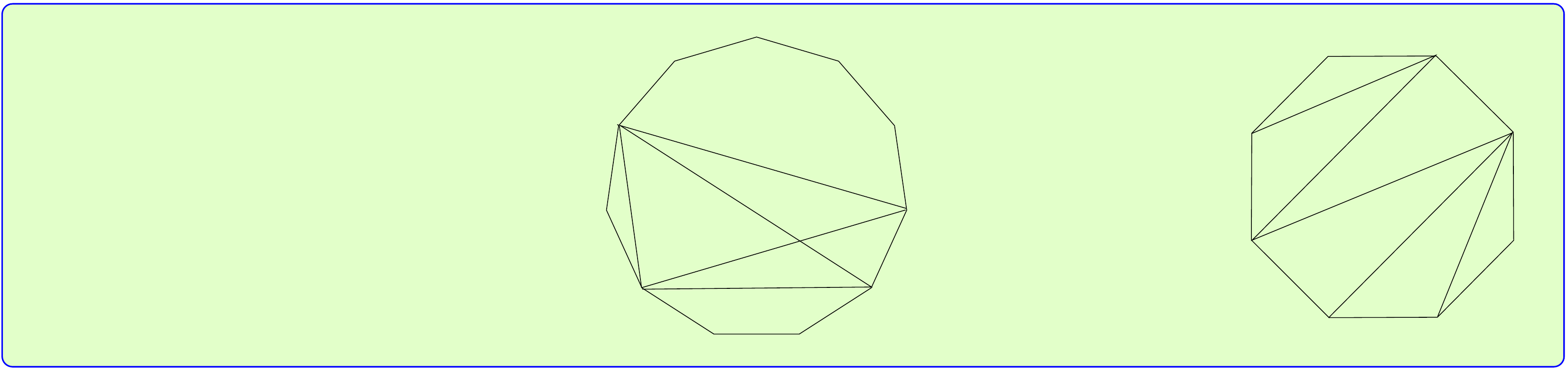}   
 \put(-430,65){\small $\begin{pmatrix}a_{11}\ a_{12}\dots a_{ii}\dots a_{1j}\dots a_{1n}\\a_{21}\ a_{22}\dots a_{2i} \dots a_{2j}\dots a_{2n} \end{pmatrix}  $}
\put(-420,30){\small ${\color{blue}p_{ij}}=det\begin{pmatrix}a_{1i} &a_{1j}\\a_{2i} & a_{2j}  \end{pmatrix}  $}    
\put(-180,10){\small  \color{blue} $p_{ik}p_{jl}=p_{ij}p_{kl}+p_{jk}p_{ki}$}
\put(-257,85){\small $n$}
\put(-221,93){\small $1$}
\put(-201,85){\small $2$}
\put(-186,67){\small $...$}
\put(-182,44){\small $i$}
\put(-192,20){\small $j$}
\put(-266,20){\small $k$}
\put(-269,67){\small $l$}
\put(-188,32){\small \color{blue} $p_{ij}$}
\put(-230,16){\small \color{blue} $p_{jk}$}
\put(-260,45){\small \color{blue} $p_{kl}$}
\put(-230,60){\small \color{blue} $p_{il}$}
\put(-243,32){\small \color{blue} $p_{ik}$}
\put(-227,45){\small \color{blue} $p_{jl}$}

\caption{Pl\"ucker coordinates, Pl\"ucker relation and a triangulation. }
\label{pluck}
\end{center}
\end{figure}
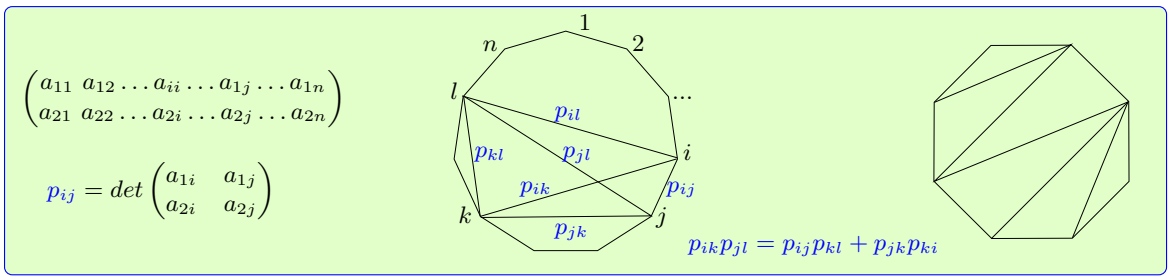

Now, let $1,\dots, n$ be the vertices of a regular $n$-gon (see Fig.~\ref{pluck}). Assign $p_{ij}$ to the diagonal of the $n$-gon connecting the vertices $i$ and $j$. Then Pl\"ucker relation will look identical to the Ptolemy relation above.

How many Pl\"ucker coordinates are needed in order  to  specify a point in  $Gr_{2,n}$? Pl\"ucker relations say that one does not need all of $p_{ij}$. More precisely, given a triangulation of the $n$-gon (see Fig.~\ref{pluck}, right) the Pl\"ucker coordinates associated to the diagonals of the triangulation and the sides of the polygon are sufficient to compute all other Pl\"ucker coordinates by applying Pl\"ucker relations several times.

\section{Conway-Coxeter Frieze Patterns}

\begin{figure}[!h]
  \begin{center}
    \includegraphics[width=0.99\linewidth]{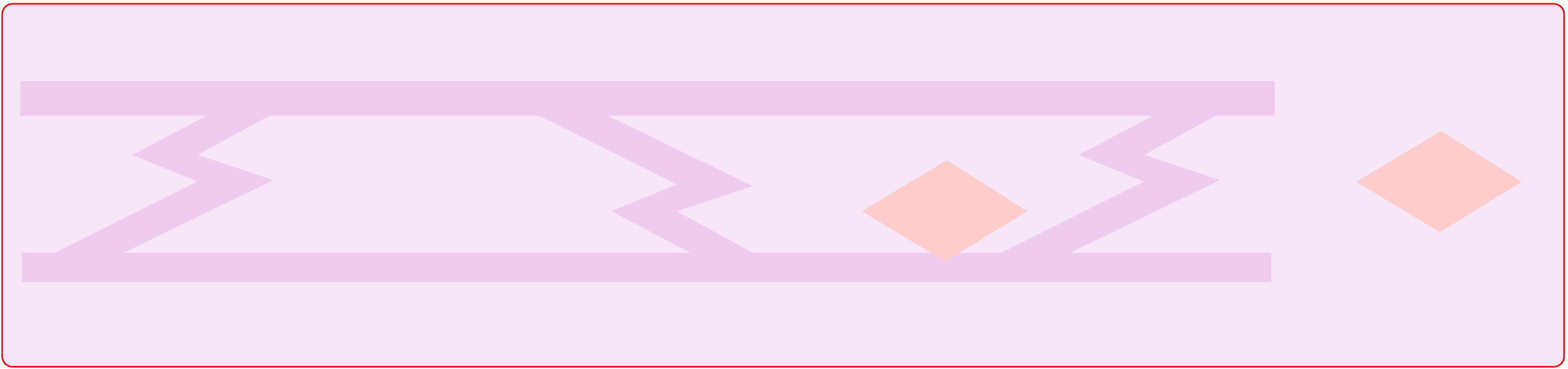}
    \put(-435,50){\scriptsize\begin{tabular}{ccccccccccccccccccccc}
$\dots$\!&&1&&1&&1&&1&&1&&1&&1&&1&&1&&$\dots$ \\
&2&&{\color{red}1}&&3&&4&&{\color{red}1}&&2&&2&&3&&2&&{\color{red}1}& \\
$\dots$\!&&{\color{red}1}&&2&&11&&3&&{\color{red}1}&&3&&5&&5&&{\color{red}1}&&$\dots$ \\
&2&&{\color{red}1}&&7&&8&&2&&{\color{red}1}&&7&&8&&2&&{\color{red}1}& \\
$\dots$\!&&{\color{red}1}&&3&&5&&5&&{\color{red}1}&&2&&11&&3&&{\color{red}1}&&$\dots$ \\
&{\color{red}1}&&2&&2&&3&&2&&{\color{red}1}&&3&&4&&{\color{red}1}&&2& \\
$\dots$\!&&1&&1&&1&&1&&1&&1&&1&&1&&1&&$\dots$ \\
                             \end{tabular}}
\put(-60,50){\scriptsize\begin{tabular}{ccc}
&$b$&\\
$a$&&$d$\\                          
&$c$&\\                                                     
           \end{tabular}}
\put(-55,25){\scriptsize $ad-bc=1$}                      
\caption{A frieze with a highlighted  diamond and zig-zags of 1's.}
\label{frieze}
\end{center}
\end{figure}

A {\it frieze pattern} is a grid of positive integers as in Fig.~\ref{frieze} with $n+2$ rows, where the first and last rows consist of 1's and every four numbers in every small diamond \quad  $a \ \ \  d$\put(-13,3){$b$}\put(-13,-3){$c$} \quad satisfy the following {\it diamond rule}: $ad-bc=1$. 
     
Conway and Coxeter~\cite{CC1,CC2} showed  that every frieze pattern %with finite number of rows
is periodic.
%Moreover, it always contains a zig-zag of 1's connecting the boundary rows.
Moreover, if one starts with two boundary rows of 1's and a connecting  zig-zag of 1's, one can always build a frieze pattern by reconstructing the entries one by one using the diamond rule. And the frieze constructed in this way will be always periodic, with period $n+3$, and will always consist of positive integers.

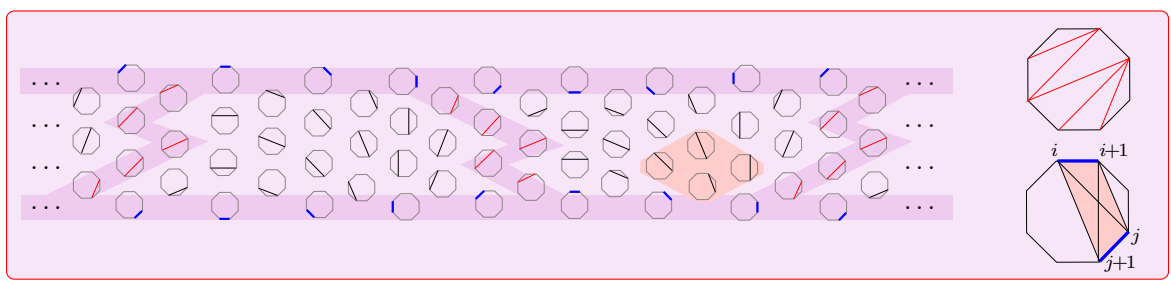
\begin{figure}[!h]
  \begin{center}
    \includegraphics[width=0.99\linewidth]{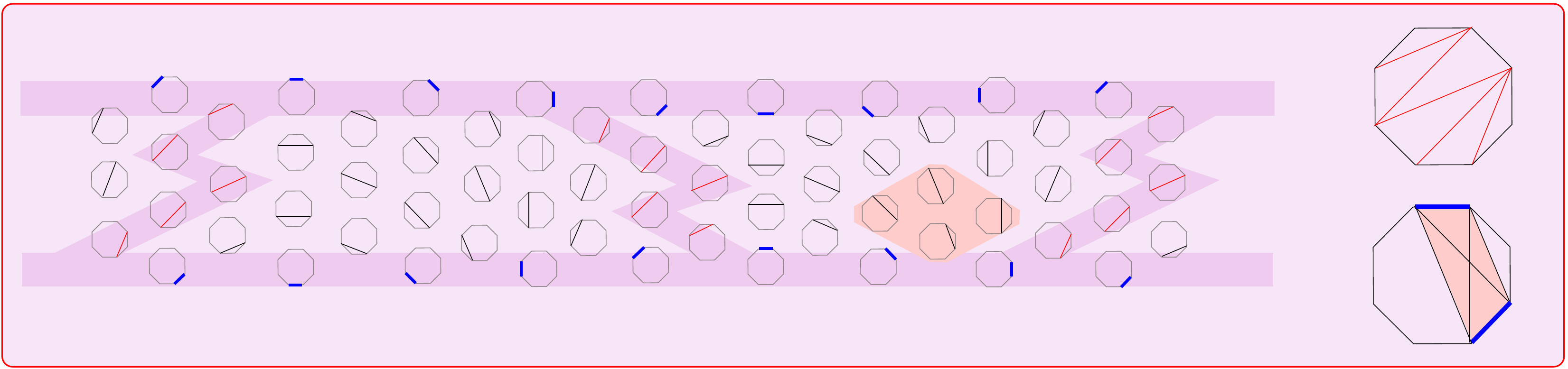}
    \put(-100,74){$\dots$}
    \put(-100,58){$\dots$}
    \put(-100,42){$\dots$}
       \put(-100,27){$\dots$} 
    \put(-427,74){$\dots$}
    \put(-427,58){$\dots$}
    \put(-427,42){$\dots$}
    \put(-427,27){$\dots$}
    \put(-15,15){\scriptsize $j$}
    \put(-25,5){\scriptsize $j\!\!+\!\!1$}
    \put(-45,47){\scriptsize $i$}
    \put(-27,47){\scriptsize $i\!\!+\!\!1$}
 \caption{Entries of the frieze and diagonals of the polygon. }
\label{frieze2}
\end{center}
\end{figure}

We will now associate the diagonals of an $(n+3)$-gon to entries of the frieze pattern (see Fig.~\ref{frieze2}) so, that the boundary edges of the polygon will correspond to  entries in the first and last rows of the frieze, shifting along the frieze to the right will correspond to a clockwise rotation of the polygon, and the adjacent entries of a diamond in the frieze correspond to diagonals with one common vertex and the other vertex shifted by one position.
A zig-zag of 1's will in this way correspond to assigning 1's to the diagonals in some triangulation of the polygon (and every triangulation will correspond to a set of entries from which the rest of the frieze can be found by the diamond rule).

Notice that the diamonds in a frieze correspond to four diagonals in the polygon. More precisely, labelling the vertices of the polygon  $1,2,3,\dots, n$ clockwise, we get diagonals of the quadrilateral $i,i+1,j,j+1$.
So, assuming that the entries of the frieze are assigned to the corresponding diagonals, we see that the diamond rule is exactly taking the shape of the Ptolemy's identity  (given that the boundary sides of the polygon are assigned 1's).

If there existed a cyclic triangulated Euclidean polygon with all sides of length 1 and all diagonals in the triangulation also of length 1,  the entries of the frieze would represent Euclidean lengths of all the diagonals of the polygon.
However, there are no  Euclidean polygons with unit sides and unit diagonals in a triangulation.
Luckily, one can overcome this by considering ideal hyperbolic polygons instead of Euclidean ones.

\section{Hyperbolic Ptolemy and triangulated surfaces}

Recall that the upper half-plane model of hyperbolic plane $\H^2$ is the set $\{z\in \C \ | \ Im(z)>0 \}$ with  metric given by $ds^2=\frac{dx^2+ dy^2}{y^2}$. In this model, geodesics are represented by half-lines and half-circles orthogonal to the real axis. Applying the map $f(z)=\frac{z-i}{z+i}$ one can transform the upper half-plane model to the Poincar\'e disc model in the unit disc. See for example~\cite{S} for more detail. 

\begin{figure}[!h]
  \begin{center}
    \includegraphics[width=0.99\linewidth]{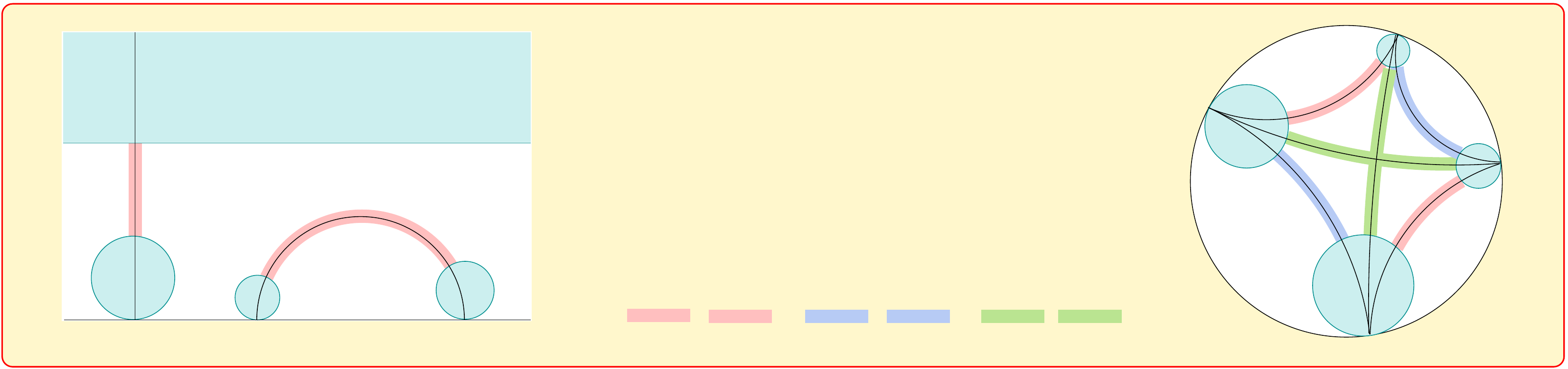}
    \put(-260,19){\scriptsize \color{blue}$\lambda_{AB}\cdot \lambda_{CD}+\lambda_{BC}\cdot \lambda_{DA}=\lambda_{AC}\cdot \lambda_{BD}$}
    \put(-230,55){{\color{blue} $\lambda_{AB}$}$=e^{l_{AB}/2}$}
    \put(-50,2){\small \color{ForestGreen} $A$}
    \put(-110,72){\small\color{ForestGreen} $B$}
    \put(-43,93){\small\color{ForestGreen} $C$}
    \put(-16,56){\small\color{ForestGreen} $D$}
    \put(-403,4){\small \color{ForestGreen} $A$}
    \put(-395,86){\small \color{ForestGreen} $B=\infty$}
    \put(-368,4){\small \color{ForestGreen} $C$}
    \put(-306,4){\small \color{ForestGreen} $D$}
    \put(-394,47){\small  $l_{AB}$}
    \put(-343,47){\small  $l_{CD}$}
    \caption{Lambda-lengths: removing infinity by horocycles} 
\label{horo}
\end{center}
\end{figure}
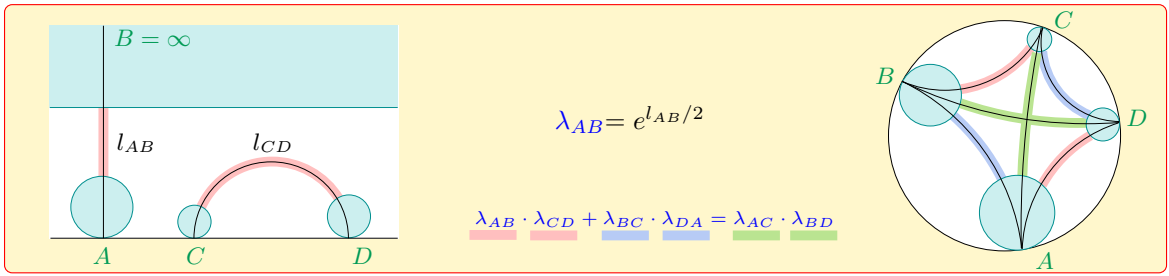

Let $A,B\in \partial \H^2$ be two points at the boundary of the hyperbolic plane. The hyperbolic distance between them is infinite, but the infinity is concentrated around $\partial \H^2$ and can be dealt with by using
  horocycles as follows (a {\it horocycle} is a limit of a circle as the centre approaches $\partial \H^2$, in the upper half-plane model of $\H^2$, a horocycle centred at $\infty$ is represented by a horizontal line, horocycles centred at other points are represented by circles tangent to  $\partial \H^2$). Choose  horocycles $h_A$ and $h_B$ centred at $A$ and $B$ (see Fig.~\ref{horo}   for the pictures in the upper half-plane  and in the disc models). Let $l_{AB}$ be the signed distance between the horocycles  $h_A$ and $h_B$ ($l_{AB}$ is zero when the horocycles are tangent and negative when the horocycles intersect each other).
Denote $\lambda_{AB}=e^{l_{AB}/2}$, the {\it lambda length} of $AB$.

An {\it ideal} polygon in $\H^2$ is a polygon with all vertices at  $\partial \H^2$. 
Given an ideal quadrilateral $ABCD$  and a choice of horocycles around each of the vertices, one can prove that the lambda lengths for $ABCD$ satisfy Ptolemy identity~\cite{Pen}:
$$ \lambda_{AB}\cdot \lambda_{CD}+\lambda_{BC}\cdot \lambda_{DA}=\lambda_{AC}\cdot \lambda_{BD}.
$$

It is easy to show that given $a,b,c>0$, one can construct (in a unique way up to isometry) an ideal hyperbolic triangle with a choice of horocycles at its vertices such that $a=\lambda_{BC}, b=\lambda_{AC}, c=\lambda_{AB}$.
Also, triangles can be attached to each other along the edges with the same lambda lengths.
Therefore, given a triangulated polygon with  positive
numbers assigned to its sides and diagonals in the triangulation, one can construct an ideal hyperbolic polygon with horocycles assigned to its vertices  such that for every diagonal the assigned number will coincide with the corresponding lambda length.
In particular, this means that every frieze pattern described above can be modelled by a hyperbolic polygon.

In a similar way, given a triangulated surface $S$   and a set of positive numbers assigned to the arcs of the triangulation, one can define a unique  hyperbolic structure (with a unique choice of horocycles at all vertices of the triangulation), so that
% the inner marked points will turn into  cusps and
the assigned numbers will coincide with the lambda lengths of the corresponding arcs.
In other words, lambda lengths of arcs in a triangulation of $S$ provide coordinates on the decorated Teichm\"uller space of $S$.

\section{Friezes from surfaces}
Let $(S,M)$ be a surface $S$ with a set of marked points $M$. We assume that it is possible to triangulate $S$ so that
the vertices of each triangle are marked points. Let $E$ be the set of arcs on $(S,M)$. One can generalise the definition of a frieze pattern from the settings of a polygon to a general marked surface in the following way.

A {\it frieze $\varphi$ from $(S,M)$} is a map $\varphi:E\to \R$ assigning a real number $\varphi(\gamma)$ to every arc $\gamma\in E$ in such a way that the Ptolemy relation holds for every quadrilateral on the surface.
A frieze  $\varphi$ is
\begin{itemize}
\item[-] {\it positive} if all numbers  $\varphi(\gamma)$ are positive;
\item[-]  {\it integer} if all numbers  $\varphi(\gamma)$ are integers;
\item[-]  {\it unitary} if there exists a triangulation $T$ of $(S,M)$ such that $\varphi(\gamma)=1$ for all $\gamma\in T$.

\end{itemize}

\medskip
\noindent
{\bf Question:} Given a marked surface $(S,M)$, is it true that every positive integer frieze  from $(S,M)$ is unitary?

When $S$ is a disc with $n$ boundary points (i.e. a polygon), Conway-Coxeter's theorem~\cite{CC1,CC2} gives a positive answer. So, one could expect this would be the case in general. However,  it was shown in~\cite{FP} that the answer is negative for a punctured polygon (friezes of type $D$). Later, a positive answer was obtained for an annulus~\cite{GS} and for a pair of pants~\cite{CGFT}.
Besides that, to our knowledge the question of unitarity of positive integer friezes from surfaces is  open.

\section{Fugue: Cluster algebras}

All instances appearing above can be described as various manifestations of cluster algebras. We will sketch the idea of cluster algebra and illustrate it by an example of cluster algebras arising from surfaces.

\subsection{Quivers, seeds, clusters...}

\noindent
{\bf Quiver mutations.} We start with a {\it quiver} $Q$, i.e. an oriented graph with finitely many vertices labelled $1,\dots,n$ and $b_{ij}\in \Z$ arrows from vertex $i$ to vertex $j$. We assume that  the quiver contains no loops ($b_{ii}=0$) and no 2-cycles (here, $p$ arrows from $i$ to $j$ are understood as $-p$ arrows from $j$ to $i$, so that $b_{ij}=-b_{ji}$).
When $b_{ij}=1$, we omit the label on the corresponding arrow of the quiver.

Next, we will pick a vertex $k$ of the quiver and will define a {\it mutation} $\mu_k$ of $Q$ taking $Q$ to $Q'=\mu_k(Q)$, which is only different from $Q$  in a small neighbourhood of the vertex $k$ and coincides with $Q$ away from $k$. The effect of  $\mu_k$ on $Q$ can be explained in two steps:

\noindent
\includegraphics[width=0.99\linewidth]{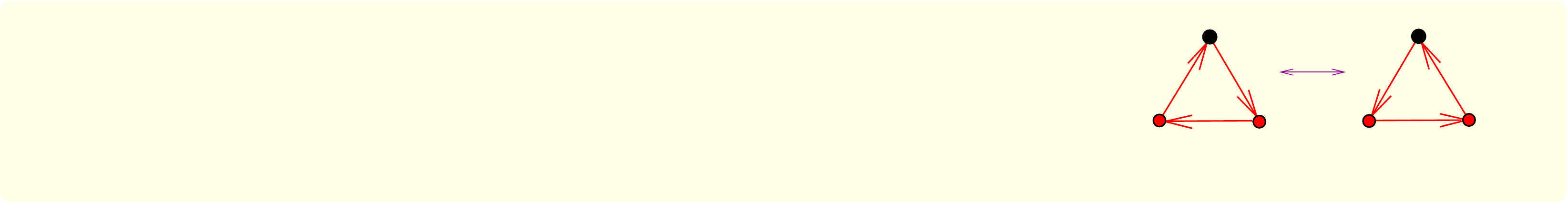}
\put(-96,47){$k$}
\put(-38,47){$k$}
\put(-75,42){\color{Plum}$\mu_k$}
\put(-92,36){\color{Blue}$q$}
\put(-114,34){\color{Blue}$p$}
\put(-100,16){\color{Blue}$r$}
\put(-33,35){\color{Blue}$q$}
\put(-55,35){\color{Blue}$p$}
\put(-45,14){\color{Blue}$r'$}
\put(-110,5){ {\color{Blue} where $r+r'=pq$}}
    \put(-420,25){\small \color{Blue} \begin{tabular}{l} (1)  reverse all arrows \\\phantom{(1)} incident to the vertex $k$;\end{tabular}}
    \put(-270,25){\small \color{Blue} \begin{tabular}{l} (2)  for every path $i\stackrel{p}\to k\stackrel{q}\to j$  \\ \phantom{(1)} with $p,q>0$ apply:\end{tabular}}

\noindent    
which corresponds to the following formula
$$b_{ij}'=\begin{cases} -b_{ij}, & \text{if $i=k$ or $j=k$} \\ b_{ij}+\frac{1}{2}(|b_{ik}|b_{kj}+b_{ik}|b_{kj}|), & \text{otherwise.} \end{cases}$$

\begin{tcolorbox}[width=\textwidth,colframe=red,colback=blue!10,boxrule=1pt] 
  \begin{ex}
    \raisebox{-20pt}{ \qquad \qquad \qquad \qquad \includegraphics[width=0.5\linewidth]{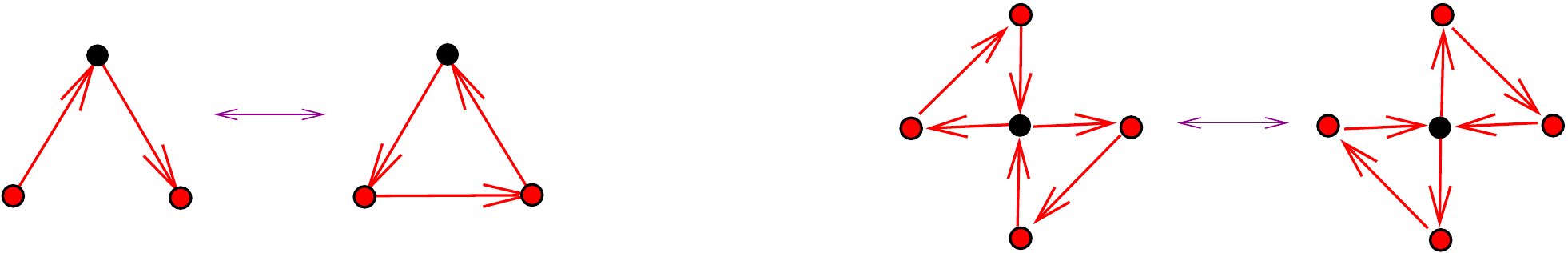}
\put(-209,1){\scriptsize $1$}
\put(-191,28){\scriptsize $2$}
\put(-180,1){\scriptsize $3$}
\put(-173,23){\small \color{Plum}$\mu_2$}
\put(-163,1){\scriptsize $1$}
\put(-145,28){\scriptsize $2$}
\put(-134,1){\scriptsize $3$}
\put(-79,0){\scriptsize $1$}
\put(-14,0){\scriptsize $1$}
\put(-93,18){\scriptsize $2$}
\put(-34,19){\scriptsize $2$}
\put(-22,19){\scriptsize $5$}
\put(-69,19){\scriptsize $5$}
\put(-57,19){\scriptsize $4$}
\put(-1,19){\scriptsize $4$}
\put(-23,30){\scriptsize $3$}
\put(-69,30){\scriptsize $3$}
\put(-48,23){\small \color{Plum}$\mu_5$}

    }
\put(-350,-16){\color{red} \bf Quiver mutation.}  
\end{ex}
\end{tcolorbox}

\medskip
\noindent
{\bf Seed mutations.}
A {\it seed} is a pair $(Q,\mathbf u)$, where $Q$ is a quiver and $\mathbf u=(u_1,\dots, u_n)$ is an $n$-tuple of algebraically independent rational functions $u_i(x_1,\dots,x_n)$ in the variables $x_i$. A function $u_i$ is associated   with the vertex $v_i$ of $Q$.
In the {\it initial seed}, we assume $u_i=x_i$.

%\begin{tcolorbox}[width=\textwidth,colback=blue!10,boxrule=0pt] 
%  \begin{ex}

%  \end{ex}  
%\end{tcolorbox}
\medskip
\noindent
A mutation $\mu_k$ in the direction $k$ will take the seed $(Q, \mathbf u)$  to  $(Q', \mathbf u')$, where $Q'=\mu_k(Q)$ is obtained by the quiver mutation described above and $\mathbf u'=(u_1',\dots, u_n')$, where $u_i'=u_i$ for all $i\ne k$ and
$$u_k'=\frac{\prod\limits_{i\to k} u_i^{b_{ik}} +\prod\limits_{j \leftarrow k} u_j^{b_{kj}}}{u_k},$$
where the products are taken over all vertices $i$ such that there is an incoming (respectively outgoing) arrow $i\to k$ (respectively
$j \leftarrow k$) with positive weight $b_{ik}>0$ (respectively, $b_{kj}>0$). 
The above relation is called {\it exchange relation}. Note that a mutation is an involution, i.e. $\mu_k(\mu_k(Q,\mathbf u))=(Q,\mathbf u)$.

\begin{tcolorbox}[width=\textwidth,colframe=blue,colback=blue!10,boxrule=1pt] 
  \begin{ex}
  \label{ex2}  
  \raisebox{-20pt}{ \qquad \qquad \qquad \qquad \includegraphics[width=0.5\linewidth]{pic/mut-ex.pdf}
    \put(-209,1){\scriptsize  \color{blue}$x_1$}
\put(-191,28){\scriptsize \color{blue} $x_2$}
\put(-180,1){\scriptsize  \color{blue}$x_3$}
\put(-173,23){\small \color{Plum}$\mu_2$}
\put(-163,1){\scriptsize  \color{blue} $x_1$}
\put(-143,26){\scriptsize  \color{blue}$x_2'=\frac{x_1+x_3}{x_2}$}
\put(-134,1){\scriptsize  \color{blue}$x_3$}
\put(-82,0){\scriptsize  \color{blue}$x_1$}
\put(-14,0){\scriptsize  \color{blue}$x_1$}
\put(-97,18){\scriptsize  \color{blue}$x_2$}
\put(-34,20){\scriptsize  \color{blue}$x_2$}
\put(-14,9){\scriptsize  \color{blue}$x_5'$}
\put(-82,10){\scriptsize  \color{blue}$x_5$}
\put(-60,19.5){\scriptsize  \color{blue}$x_4$}
\put(-1,19){\scriptsize  \color{blue}$x_4$}
\put(-27,30){\scriptsize  \color{blue}$x_3$}
\put(-69,30){\scriptsize  \color{blue}$x_3$}
\put(-48,23){\small \color{Plum}$\mu_5$}
\put(6,5){\small \color{blue}$x_5'=\frac{x_1x_3+x_2x_4}{x_5}$}
  }
\put(-350,-16){\color{blue} \bf Seed mutation.}  
\end{ex}
\end{tcolorbox}

The functions $u_i(x_1,\dots,x_n)$ obtained from  $x_1,\dots,x_n$ in the process of interated seed mutations are called
{\it cluster variables}. The collection of cluster variables contained in one seed is called a {\it cluster}.
Let  $\Q(x_1,\dots,x_n)$ be  the field of rational functions in the variables $x_i$, $i\in\{1,\dots,n\}$.
The {\it cluster algebra} $\mathcal A (Q)$ is  the $\Q$-subalgebra of $\Q(x_1,\dots,x_n)$   generated by all cluster variables.

\begin{tcolorbox}[width=\textwidth,colframe=Green, colback=blue!10,boxrule=1pt] 
  \begin{ex}
    \raisebox{-60pt}{ \qquad \qquad \qquad \qquad \qquad \includegraphics[width=0.6\linewidth]{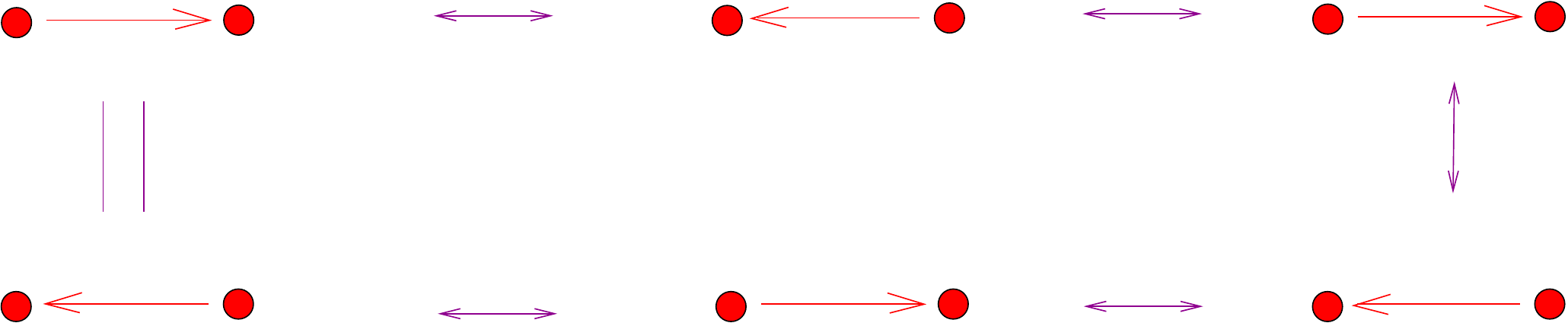}    
\put(-255,52){\color{blue} \small $x_1$}
\put(-206,52){\color{blue} \small $x_2$}
\put(-255,6){\color{blue} \small $x_2$}
\put(-206,6){\color{blue} \small $x_1$}
\put(-143,52){\color{blue} \small $x_1$}
\put(-103,58){\color{blue} \small  $\frac{1+x_2}{x_1}$}
\put(-143,6){\color{blue} \small $x_2$}
\put(-103,12){\color{blue} \small  $\frac{1+x_1}{x_2}$}
\put(-58,12){\color{blue} \small  $\frac{1+x_1+x_2}{x_1x_2}$}
\put(-7,12){\color{blue} \small  $\frac{1+x_1}{x_2}$}
\put(-58,58){\color{blue} \small  $\frac{1+x_1+x_2}{x_1x_2}$}
\put(-7,58){\color{blue} \small  $\frac{1+x_2}{x_1}$}
\put(-15,29){\color{Plum} \small  $\mu_2$}
\put(-71,53){\color{Plum} \small  $\mu_1$}
\put(-71,7){\color{Plum} \small  $\mu_1$}
\put(-171,53){\color{Plum} \small  $\mu_2$}
\put(-171,7){\color{Plum} \small  $\mu_2$}
}
\put(-410,-15){\color{Green} \bf Cluster algebra of type $A_2$:}  
\put(-410,-28){5 seeds, 5 cluster variables}  
\put(-410,-45){\small \color{Plum} (the seed after 5 mutations }  
\put(-410,-55){\small \color{Plum} coincides with the initial seed }  
\put(-410,-65){\small \color{Plum} up to a permutation of the vertices)}  
\end{ex}
\end{tcolorbox}

\vspace{15pt}
\noindent
{\bf Laurent Phenomenon (\cite{FZ1}):} In a cluster algebra, every cluster variable  $u_i(x_1,\dots,x_n)$ is a Laurent polynomial in $(x_1,\dots,x_n)$, i.e. $u=\frac{P(x_1,\dots,x_n)}{Q(x_1,\dots, x_n)}$, where $Q(x_1,\dots,x_n)=x_1^{d_1}\dots x_n^{d_n}$ is a monomial.

\vspace{10pt}

\subsection{Examples of cluster algebras}

\begin{tcolorbox}[width=\textwidth,colback=blue!10,boxrule=0pt] 
  \begin{ex}[Cluster algebra from triangulated surfaces]
 \label{ex4}   
%\begin{tcolorbox}[width=\textwidth,colback=blue!10,boxrule=1pt]   
    Given a triangulated surface, one can construct a quiver as shown below.
  Each vertex $i$ of the quiver here corresponds to some arc $e_i$ of the triangulation.
  It is easy to check that a mutation $\mu_i$ of the quiver corresponds to a flip of the corresponding edge $e_i$ in the triangulation.
  \begin{center}
    \includegraphics[width=0.99\linewidth]{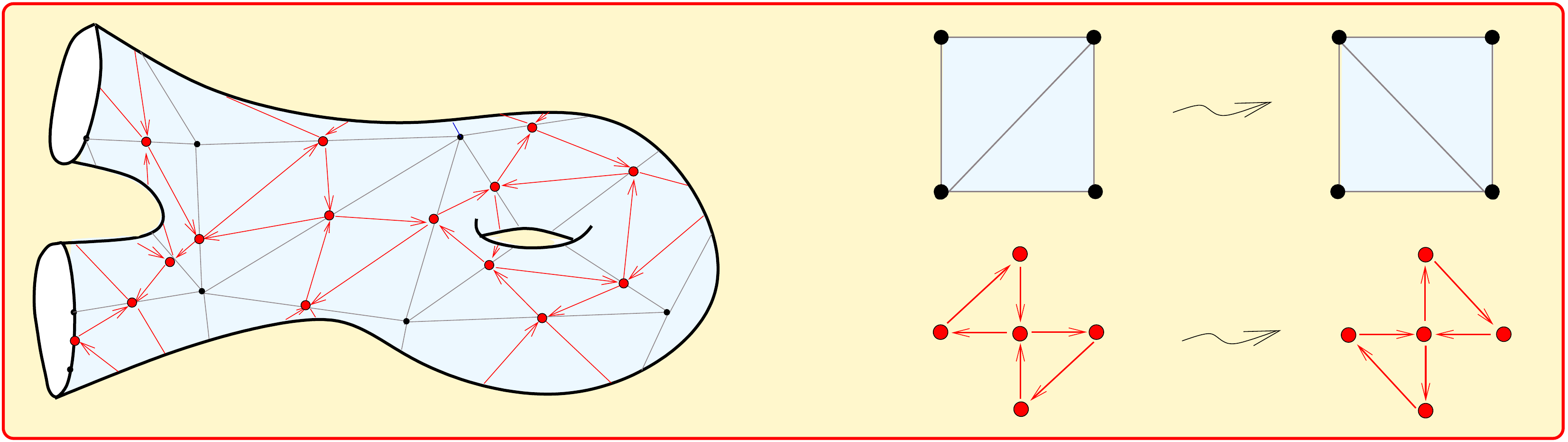}
    \put(-152,87){$e_k$}
    \put(-152,87){$e_k$}
    \put(-102,95){flip $f_k$}
   \put(-102,43){mutation}
   \put(-90,34){$\mu_k$}
   \put(-140,30){$k$}
\end{center}
%\end{tcolorbox}

\begin{tcolorbox}[width=0.99\textwidth,colframe=blue, colback=blue!10,boxrule=1pt] 
%\begin{ex}
  For a triangulated surface, the part of the quiver $Q$ adjacent to $k$ (and thus defining the mutation $\mu_k$) is as in Example~\ref{ex2} , and the {\color{blue} exchange relation} turns into
$${\color{blue}   u_k'=\frac{u_iu_j+u_lu_m}{u_k}}, $$
  which is identical to   Ptolemy (or to Pl\"ucker) relation.  
\end{tcolorbox}

\begin{tcolorbox}[width=0.99\textwidth,colframe=Green, colback=blue!10,boxrule=1pt] 
  Given a triangulated surface and initial values $x_1,\dots,x_n$ assigned to the arcs of the triangulation, one can construct a unique hyperbolic surface with a choice of horocycles around vertices, so that $x_i$ will coincide with the lambda lengths of the corresponding arcs. A sequence of iterated mutations of the initial seed will then correspond to applying a sequence of flips  to a triangulation,
cluster variables as functions of $\{x_i\}$   will coincide with the lambda lengths of the corresponding arcs, clusters will correspond to collections of lambda lengths of arcs in one triangulation. See~\cite{FST} and~\cite{FT}.
\end{tcolorbox}

\end{ex}
\end{tcolorbox}

%\vspace{10pt}

%\begin{tcolorbox}[width=\textwidth,colback=SpringGreen!30,boxrule=0pt] 
\begin{tcolorbox}[width=\textwidth,colback=LimeGreen!30,boxrule=0pt] 
  \begin{ex}
    (Cluster algebra from Grassmannian $G_{2,n}$). Given an  $n$-gon $P$ with a triangulation $T$, one can retell the same story as in Example~\ref{ex4}.
%\begin{tcolorbox}[width=0.98\textwidth,colframe=blue,colback=SpringGreen!30,boxrule=0.5pt]
A  {\color{red}quiver $Q$} is built from the triangulation: the vertices of $Q$ correspond to the diagonals and sides of $P$.
%\end{tcolorbox}
  A {\color{blue} cluster variable} associated to the diagonal connecting vertices $i$ and $j$ is now identified with the Pl\"ucker coordinate {\color{blue} $p_{ij}$}. The {\color{blue} exchange relation} takes the form {\color{blue}$p_{jl}=\frac{p_{ij}p_{kl}+p_{jk}p_{ki}}{p_{ik}}$} and can be identified with the Pl\"ucker relation.
We obtain a {\color{Green} cluster algebra} in the ring of rational functions on  $G_{2,n}$, where 
 the  {\color{Green} clusters} correspond to the collections of Pl\"ucker coordinates associated with diagonals in one triangulation. 
    
  \end{ex}  
\end{tcolorbox}

\begin{tcolorbox}[width=\textwidth,colback=Plum!10,boxrule=0pt] 
  \begin{ex} (Friezes and cluster category)
    Cluster categories were introduced in~\cite{BMRRT} %[CCS06]
    to give a categorical interpretation of cluster algebras.
    For a cluster algebra associated with a triangulated polygon (with a good triangulation),   a cluster category can be represented by a frieze where each entry is replaced  by a finitely generated module over $kQ$ (here $Q$ is the quiver corresponding to the triangulation and $k$ is an algebraically closed field). %The modules in one diamond are connected by morphisms.
    There is a tight connection between a cluster algebra and the corresponding cluster category, in particular, many results about cluster algebras were established using cluster categories. See~\cite{MG} for a review about friezes and~\cite{Pr} for a bridge to cluster categories.

\end{ex}
\end{tcolorbox}

%Cluster variables have many surprising properties, in particular:

%\begin{itemize}
%\item {\it Laurent phenomenon} (conjectured and proved in~\cite{FZ1},~\cite{FZ2}): In a cluster algebra, all cluster variables are Laurent polynomials in $(x_1,\dots,x_n)$, i.e. a cluster variable has a shape $\frac{P(x_1,\dots,x_n)}{x_1^{i_1}\dots x_n^{i_n}}$. 
%
%\item {\it Positivity} (conjectured in ~\cite{FZ1}, proved in~{LS,GHKH}): The numerator $P(x_1,\dots,x_n)$ above is a polynomial with positive coefficients.
%\end{itemize}
  
%Notice that these are global properties which hold for all seeds in a cluster algebra. At the same time, the defining conditions on cluster algebra are local, i.e. to show that there is a cluster algebra corresponding to some iterative process, one only needs to check that one iteration may be considered as mutation.

\subsection{Conclusion}
Ptolemy's theorem is deeply tied with the recent theory of cluster algebras.
Since the introduction of cluster algebras by Fomin and Zelevinsky,  cluster algebras found connections and applications to  a large range of domains in mathematics and mathematical physics. This includes (but is not exhausted by!) the following:

\vspace{5pt}
  \begin{center}
    \includegraphics[width=0.99\linewidth]{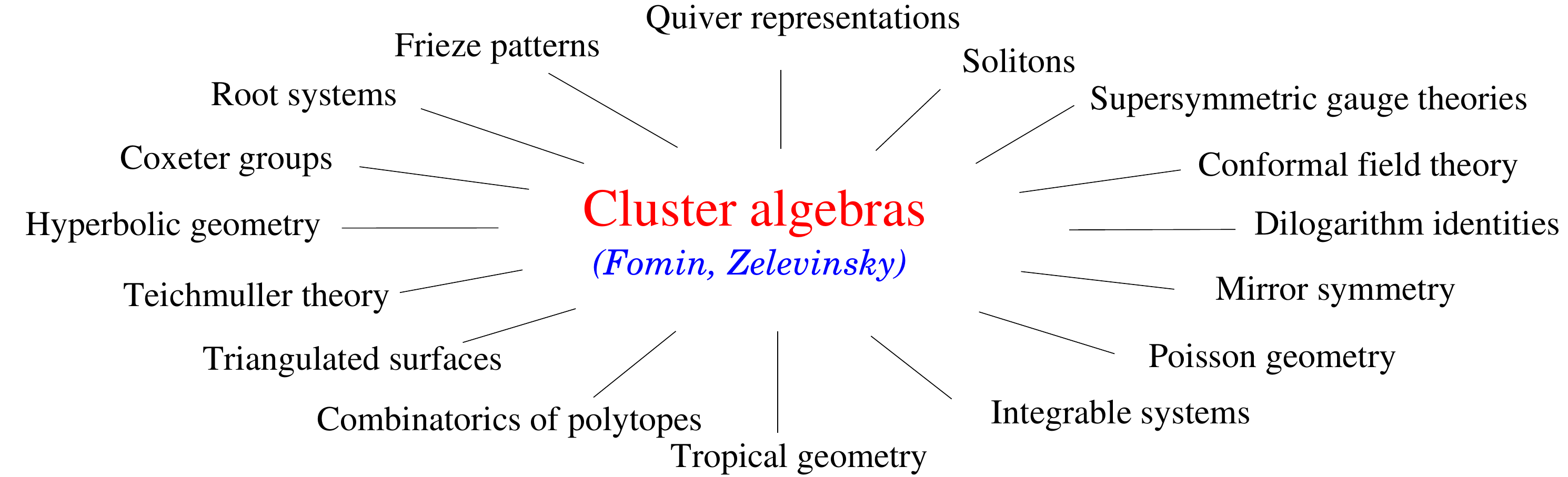}
\end{center}
\vspace{10pt}

\noindent
Links related to cluster algebras are collected on \href{http://www.math.lsa.umich.edu/~fomin/cluster.html}{Cluster Algebras Portal}~\cite{F} maintained by Sergey Fomin. One can start reading with~\cite{W},~\cite{M} and~\cite{FR}.

\begin{rem}
Another  train of results related to the Ptolemy relation arises from the Ptolemy inequality stating that for any quadrilateral $ABCD$ in the Euclidean plane,  $AC\cdot BD \le AB\cdot CD+BC\cdot AD$. Metric spaces with this property are called {\it Ptolemy spaces}.
A metric space is CAT(0) if and only if it is Ptolemy and
Busemann convex, see~\cite{FLS}.

\end{rem}

\section{Back to Ptolemy: proof without words}
There are \href{https://ckrao.wordpress.com/2015/05/24/a-collection-of-proofs-of-ptolemys-theorem/}{numerous proofs} of Ptolemy's theorem and Ptolemy's inequality based on various ideas (similarity, inversion, triangle inequality, etc).
To conclude, we  provide a ``proof without words''   (appeared in~\cite{DH} and popularised by the \href{https://www.cut-the-knot.org/proofs/PtolemyTheoremPWW.shtml}{cut-the-knot} portal).  

\begin{figure}[!h]
  \begin{center}
\includegraphics[width=0.99\linewidth]{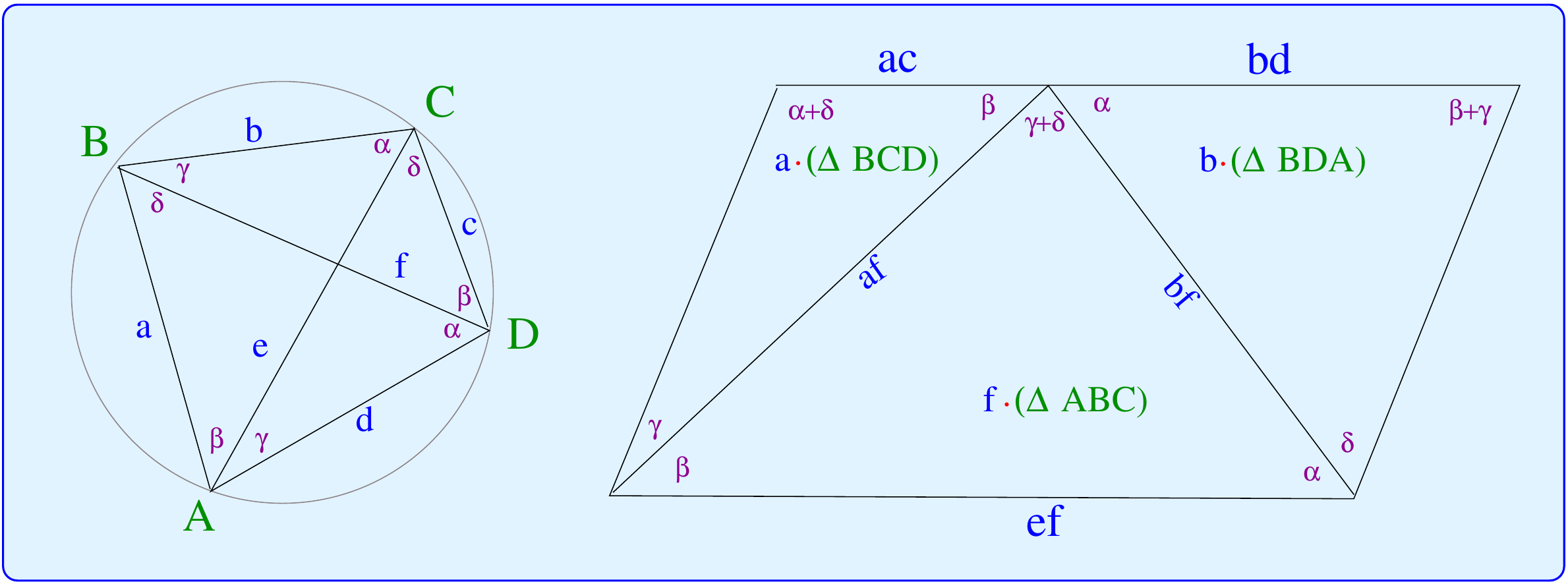}       
 \caption{Ptolemy's theorem:  {\color{blue} \bf ef=ac+bd}}
\label{ptolemy-pf}
\end{center}
\end{figure}
  
\pagebreak
\clearpage

\noindent
{\bf Acknowledgement. }  The author is greateful to Sergei Tabachnikov for encouragement; to Pavel Tumarkin, Gil Bor and Evgeny Smirnov for helpful discussions  and for useful suggestions on the early version of this note;    and to \Ilke \Canakci \  for her patience in explaining connections between clusters, freizes and cluster categories.

\end{document}